\documentclass[11pt,a4paper,leqno,twoside, headinclude]{amsart}

\usepackage[tracking=true]{microtype}
\DeclareMicrotypeSet*[tracking]{my}%
  { font = */*/*/sc/* }%
\SetTracking{ encoding = *, shape = sc }{ 45 }

\title{On quasi-isometric nilpotent Lie groups}
\def\titl{On quasi-isometric nilpotent Lie groups}
\def\auth{Manuel Amann}
\date{October 12th, 2017}

\thanks{}

\author{\auth}

\usepackage[english]{babel}

\usepackage[colorlinks,pdftex, plainpages=false]{hyperref}
\hypersetup{pdftitle=\titl, pdfauthor=\auth, pdftoolbar=false,
plainpages=false, hyperindex=true, pdfdisplaydoctitle=true}

\hypersetup{colorlinks,%
citecolor=black,%
filecolor=black,%
linkcolor=black,%
urlcolor=black,%
pdftex}

\usepackage{color}

\usepackage{amsmath, amssymb, amscd, amsthm}
\usepackage{stmaryrd}
\usepackage{latexsym}
\usepackage[all]{xy}
\usepackage{pb-diagram}
\usepackage{rotating}
\usepackage{multicol}
\usepackage{lscape}
\usepackage{wasysym}
\usepackage{longtable}
\usepackage{enumerate}

\xyoption{all}

\newtheorem{theo}{Theorem}[section]
\newtheorem{main}{Theorem}

\newtheorem*{main*}{Theorem}

\newtheorem*{mainprop*}{Proposition}

\newtheorem{mainconj}{Conjecture}

\newtheorem{defi2}[theo]{Definition}
\newtheorem*{defi2*}{Definition}

\newenvironment{defi*}{\begin{defi2*}\normalfont}{\end{defi2*}}

\newenvironment{defin*}[1]{\begin{defi2*}[#1]\normalfont}{\end{defi2*}}
\newtheorem*{rem2*}{Remark}
\newenvironment{rem*}{\begin{rem2*}\normalfont}{\hfill$\boxbox$\end{rem2*}}
\newtheorem{rem2}[theo]{Remark}

\newtheorem{lemma}[theo]{Lemma}

\newtheorem*{cor*}{Corollary}

\newtheorem{conj}[theo]{Conjecture}

\newtheorem*{conj*}{Conjecture}
\newtheorem*{theo*}{Theorem}
\newtheorem*{ques*}{Question}
\newtheorem*{mi2}{Main Idea}

\newtheorem{ex2}[theo]{Example}

\newtheorem{exer2}[theo]{Exercise}

\newtheorem{alg2}[theo]{Algorithm}

\newcommand{\nn}{{\mathbb{N}}}                                     
\newcommand{\rr}{{\mathbb{R}}}                                     
\newcommand{\dif} {{\operatorname{d}}}                             
\newcommand{\In} {{\,\subseteq\,}}                                 
\newcommand{\co}{\colon\thinspace}                                 

\newcommand{\comment}[1]{}                                         
\newcommand{\hto}[1]{\overset{#1}{\hookrightarrow}}                

\newcommand{\case}[1]{\textbf{Case #1.}}                           
\newcommand{\step}[1]{\textbf{Step #1.}}                           
\newcommand{\ack}{\noindent\textbf{Acknowledgements. }}            
\newcommand{\str}{\noindent\textbf{Structure of the article. }}    

\newenvironment{prf}{\begin{proof}[\textsc{Proof}]} {\end{proof}}     

\begin{document}

\maketitle \thispagestyle{empty}


\begin{abstract}
In this article we provide evidence for a well-known conjecture which states that quasi-isometric simply-connected nilpotent Lie groups are isomorphic. We do so by constructing new examples which are rigid in the sense that whenever they are quasi-isometric to any other simply-connected nilpotent Lie group, the groups are actually isomorphic.
\end{abstract}


\section*{Introduction}

This article centers around the following
\begin{conj}\label{conj01}
Quasi-isometric simply-connected nilpotent Lie groups are isomorphic.
\end{conj}

We slightly rephrase the conjecture. We call a simply-connected nilpotent Lie group $G$ \emph{rigid} if any simply-connected nilpotent Lie group $G'$ quasi-isometric to $G$
is actually isomorphic to $G$.
There was a proof of Conjecture \ref{conj01} announced in the first version and withdrawn in the second version of the article \cite{KD15}. It is the goal of this article to construct new examples of rigid nilpotent Lie groups.

By a result of Pansu it is known that the associated Carnot graded nilpotent Lie algebras are isomorphic, by a result of Sauer quasi-isometric nilpotent Lie groups yield Lie algebras with isomorphic real cohomology algebras. Hence, in order to determine rigidity of a simply-connected nilpotent Lie group, it suffices to show that whenever there is a second group with isomorphic associated Carnot graded Lie algebra and isomorphic real cohomology algebra, then both Lie algebras are actually already isomorphic. We shall do so by passing to a dual description in real homotopy theory.

We provide the following new examples---drawing on the 1:1-correspondence between simply-connected Lie groups and Lie algebras. In particular, the second and fourth examples are not Carnot graded.
\begin{main}\label{theoB}
Let $k\geq 1$ for (1) and (3), suppose that $k\in \{2,3\}$ for (2). The nilpotent Lie groups belonging to the following nilpotent Lie algebras are rigid. The first and the third one is Carnot graded, the second one and the fourth one is not.
\begin{enumerate}
\item
\begin{align*}
\mathfrak{g}=\langle x_1,\ldots, x_{2k}, n_1,\ldots, n_{2k-1},m\rangle
\end{align*}
with Lie brackets trivial except for
\begin{align*}
[x_i,x_{i+1}]&=-n_i, \\
[x_i, n_i]&=-m \qquad \textrm{for all } 1\leq i\leq 2k-1
\end{align*}
\item
\begin{align*}
\mathfrak{g}=\langle x_1,\ldots, x_{2k}, x_{2k+1}, n_1,\ldots, n_{2k-1},m\rangle
\end{align*}
with Lie brackets trivial except for
\begin{align*}
[x_i,x_{i+1}]&=-n_i,\\
[x_i, n_i]&=-m\\
[x_{2k},x_{2k+1}]&=-m
 \qquad \textrm{for all } 1\leq i\leq 2k-1
\end{align*}
\item
The Lie algebra $\mathfrak{g}$ is free to degree $k$, but in its latest ``free degree'' $k+1$ it consists only of a subspace of the corresponding free nilpotent algebra. More precisely, there are no relations on the Lie brackets up to degree $k$, in degree $k+1$ there are potentially relations, in degrees above $k+1$ all Lie brackets vanish.
\item
\begin{align*}
\mathfrak{g}=\langle x_1,x_2,x_3, x_4, , n_1,\ldots, n_{5},m\rangle
\end{align*}
with Lie brackets trivial except for
\begin{align*}
[x_1,x_2]&=-n_1\\
[x_1,x_3]&=-n_2\\
[x_2,x_3]&=-n_3\\
[x_1,x_4]&=-n_4\\
[x_2,x_4]&=-n_5\\
[x_1,n_1]=[x_1,n_2]=[x_2,n_3]=[x_2,x_4]&=-m
\end{align*}
\end{enumerate}
\end{main}
It is clear, and we leave details to the interested reader, that Example (3) can be used to produce a rigid example in each dimension from a certain (small) dimension onward. As we were informed by Yves Cornulier there is a rigid example known in dimension 5 (see $L_{5,5}$ on  \cite[p.~47]{Cor16a}).

\vspace{5mm}

This result will be proved in four separate theorems dealing with the four different Lie algebras, namely Theorems \ref{theo01}, \ref{theo02}, \ref{theo03}, and \ref{theo04}.

The method of this proof is to first transcribe the situation to the realm of minimal Sullivan models. We then combine the results of Sauer and Pansu to carefully analyse our examples paralleled by an investigation of the minimal Sullivan model of a potentially quasi-isometric Lie group.

The realisation as a nilpotent Lie group the Lie algebra of which is rationally equivalent to the Lie algebra we start with is provided by \cite[Section 5]{CP00}. Thus it will suffice to investigate the corresponding Sullivan models in Section \ref{sec02}.

Finally, in Section \ref{sec03} (see Theorem \ref{theo07}) we provide a probably well-known example of two nilpotent Lie algebras which have isomorphic associated Carnot graded Lie algebras and isomorphic cohomology, but which are not isomorphic themselves. This shows that a general confirmation of the conjecture just from the results by Pansu and Sauer is a priori not possible.

\vspace{5mm}

From the point of view of rational homotopy theory we aim to find models of nilpotent Lie algebras which are so to say ``intrinsically formal'' within the class of nilpotent Lie algebras with the same associated Carnot graded algebra. That is, we want to find nilpotent Lie algebras already determined by their associated Carnot graded algebra and their cohomology type. There is a related and converse question: Yves F\'elix asked whether a rational cohomology type could be realised by non-isomorphic elliptic rational homotopy types. There are indeed examples of this form (see \cite{NSY03}, \cite{NSY04}). We reduced Conjecture \ref{conj01} to the following interesting and somewhat converse
\begin{ques*}
How many nilpotent Lie algebras are there sharing one associated Carnot graded type, yet possessing pairwise distinct cohomology rings?
\end{ques*}
Our examples at least yield some more instances upon which the corresponding decision question can be answered in the affirmative. As we were informed by Yves Cornulier for example in dimension 7 there are 4 families sharing a common associated graded algebra.

\vspace{3mm}

\str In Section \ref{sec01} we review and provide necessary tools for the proof of the main theorem. We also transcribe nilpotent Lie algebras to the language of minimal Sullivan models. In Section \ref{sec02} we prove Theorem \ref{theoB}, before, in Section \ref{sec03} we show that the methods used in this proof cannot solve the whole conjecture without further input.

\vspace{3mm}

\ack The author is very grateful to Roman Sauer for both introducing him to the conjecture and for several helpful discussions as well as valuable feedback in the different stages of the project. The author also thanks Yves Cornulier for helpful remarks on a previous version of the article.

The author is funded by a Heisenberg stipend of the German Research Foundation.


\section{Prerequisites}\label{sec01}

The results draw on two properties known to be implied by quasi-isometry (see \cite{Pan89},\cite{Sau06}).
\begin{theo}[Sauer]\label{theo04}
Quasi-isometric finitely generated nilpotent groups possess Malcev completions with isomorphic real graded cohomology algebras.
\end{theo}
Recall that by a theorem of Malcev any finitely generated torsion-free nilpotent group $\Gamma$ is discretely and cocompactly embedded in a unique simply-connected nilpotent Lie group G, the so-called \emph{(real) Malcev completion} of
$\Gamma$. Thus $\Gamma$ has an associated real Lie algebra
$\mathfrak{g}$ and the real cohomology algebras of $\mathfrak{g}$ and $\Gamma$ are isomorphic. For simply-connected nilpotent Lie groups the real Malcev completion agrees with the group itself.
\begin{theo}[Pansu]\label{theo05}
The associated Carnot graded real Lie algebras of quasi-isometric finitely generated nilpotent groups are isomorphic.
\end{theo}

Let us now describe how the minimal Sullivan model together with a lower grading arises from a real nilpotent Lie algebra $\mathfrak{g}$.

A \emph{Carnot grading} of $\mathfrak{g}$ is a decomposition $\mathfrak{g}=\bigoplus_{i\in \nn_0} \mathfrak{g}_i$ of $\rr$-submodules $\mathfrak{g}_i$ such that $[\mathfrak{g}_i,\mathfrak{g}_j]\In \mathfrak{g}_{i+j}$ for all $i,j\in \nn_0$. Further we require that $\mathfrak{g}_0$ generates $\mathfrak{g}$.

Equivalently (see \cite{Cor16}), one may describe this property as follows. Consider the lower central series of the Lie algebra $\mathfrak{g}$.
\begin{align*}
\mathfrak{g}^{(0)}&=\mathfrak{g}\\
\mathfrak{g}^{(1)}&=[\mathfrak{g},\mathfrak{g}]\\
\mathfrak{g}^{(2)}&=[[\mathfrak{g},\mathfrak{g}],\mathfrak{g}]\\
&\vdots \\
\mathfrak{g}^{(i)}&=[\mathfrak{g}^{(i-1)},\mathfrak{g}]\\
&\vdots
\end{align*}
The \emph{associated (Carnot-)graded Lie algebra $\operatorname{Car}(\mathfrak{g})$} is defined as
\begin{align*}
\operatorname{Car}(\mathfrak{g}):=\bigoplus_{i\geq 0} \mathfrak{g}^{(i)}/\mathfrak{g}^{(i+1)}
\end{align*}
with the induced Lie bracket and grading. Now, the Lie algebra $\mathfrak{g}$ is Carnot graded if and only if $\mathfrak{g}\cong \operatorname{Car}(\mathfrak{g})$.

Let us now inductively construct the real vector space underlying the real minimal Sullivan model of $\mathfrak{g}$. Suppose that the Lie algebra is $n$-step nilpotent, i.e.~$\mathfrak{g}^{(n-1)}\neq 0$ and $\mathfrak{g}^{(n)}=0$. Set $V_{n-1}:=\mathfrak{g}^{(n-1)}$.
Consider the short exact sequence
\begin{align*}
0\to \mathfrak{g}^{(i)} \to \mathfrak{g}^{(i-1)} \to \mathfrak{g}^{(i-1)}/\mathfrak{g}^{(i)}\to 0
\end{align*}
Suppose we have constructed a basis of $\mathfrak{g}^{(n)}$, then via complementary elements $v_j$ we may extend this basis to one of $\mathfrak{g}^{(n-1)}$. Denote by $V_i$ the space spanned by these complementary $v_j$. It follows that
\begin{align*}
V:=\bigoplus_{0\leq i\leq n-1}V_i\cong \bigoplus_{0\leq i\leq n-1} \mathfrak{g}^{(i)}/\mathfrak{g}^{(i+1)}\cong \mathfrak{g}
\end{align*}
(Usually, one would construct $V$ in a way such that it is isomorphic to $\mathfrak{g}^*$; however, we identify duals here.) In particular, we have that $V_0$ is isomorphic to the abelianisation of $\mathfrak{g}$.

When addressing the vector space $V=\bigoplus_{0\leq i\leq n-1} V_i$ we refer to the degree in the Carnot grading as the \emph{lower degree}. We extend this multiplicatively/additively to $\Lambda V$.

Denote by $\langle X_i\rangle_i$ a basis of $\operatorname{Car}(\mathfrak{g})=\bigoplus_{i\geq 0} \mathfrak{g}^{(i)}/\mathfrak{g}^{(i+1)}$ subject to the direct sum decomposition and in bijection with the basis $\langle v_i\rangle_i$ we constructed for $V$.

Note that by the short exact sequence above we can consider the basis $\langle X_i\rangle$ at the same time as a basis of $\mathfrak{g}$ and of $\operatorname{Car}(\mathfrak{g})$.

The differential $\dif$ on $V$ is defined by
\begin{align*}
\dif v_i&=\sum - c_{l,k}^i v_l v_k \intertext{where}
[X_l,X_k]&=\sum c_{l,k}^i X_i
\end{align*}
The Jacobi identity is equivalent to $\dif^2=0$.

Note that the differential has the property that it reduces lower degree, i.e.~
\begin{align*}
\dif V_i\In \Lambda V_{<i}
\end{align*}
This is equivalent to the fact that the Lie bracket increases the Carnot grading.

We call a minimal Sullivan model with this property and with $V=V_0\oplus \ldots \oplus V_{n-1}$ an \emph{$n$-stage} model.

We now deduce the minimal model $(\Lambda V,\bar \dif)$ of the associated graded algebra $\operatorname{Car}(\mathfrak{g})$.

Clearly, the associated graded algebra is indeed Carnot graded. It is clearly generated by $\mathfrak{g}^{(0)}$ and the induced Lie bracket $\overline{[\cdot, \cdot]}$ satisfies
\begin{align*}
\overline{[X_i,X_j]}\in \mathfrak{g}^{(i+j+1)}/\mathfrak{g}^{(i+j+2)}=(\operatorname{Car}(\mathfrak{g}))^{(i+j+1)}
\end{align*}
for $X_i\in \mathfrak{g}^{(i)}/\mathfrak{g}^{(i+1)}=(\operatorname{Car}(\mathfrak{g}))^{(i)}$, $X_j\in \mathfrak{g}^{(j)}/\mathfrak{g}^{(j+1)}=(\operatorname{Car}(\mathfrak{g}))^{(j)}$, as follows from an iterative application of the Jacobi identity.

For $v_i\in V_i$ the differential has the form
\begin{align*}
\bar \dif v_i&=\sum - c_{l,k}^i v_l v_k
\end{align*}
where $v_l, v_k\in V_{<i}$. This is equivalent to $[X_l,X_k]=\sum c_{l,k}^j X_j$ with the $X_j$ of Carnot degree at least $i$, and $i$ is the sum of the Carnot degrees of $X_l$ and $X_k$ plus $1$, i.e.~$i=l+k+1$. Consequently, in $\operatorname{Car}(\mathfrak{g})$ we project to those of the $X_j$ of Carnot degree exactly $i$. Again, conversely, this implies that the differential $\bar \dif$ on $(\Lambda V,\bar \dif)$ satisfies
\begin{align*}
\bar \dif \co V_i\mapsto \bigoplus_{k+l=i-1} V_k\cdot V_l \In \Lambda^2 V_{<i}
\end{align*}
defined by
\begin{align*}
\dif-\bar \dif\co V_i\mapsto \bigoplus_{k+l<i-1} V_k\cdot V_l \In \Lambda^2 V_{<i-1}
\end{align*}
This means, in particular, that $\bar \dif$ decreases the lower degree by exactly one in the sense that $\operatorname{lowdeg} \dif v=\operatorname{lowdeg} v -1$.
The model of the associated graded algebra is obviously minimal.

\vspace{5mm}

We may identify the real homotopy groups with $V$.
\begin{lemma}\label{lemma02}
If two real nilpotent Lie algebras have isomorphic associated graded Lie algebras, then they have isomorphic real homotopy groups, and their respective minimal models satisfy $(\Lambda (V_0\oplus V_1), \dif)\cong (\Lambda (W_0\oplus W_1),\dif)$.
\end{lemma}
\begin{prf}
Since, for a minimal model $(\Lambda V,\dif)$, the observation above reveals $(\Lambda V,\bar \dif)$ as a minimal model of the associated graded algebra, and since both algebras are free over $V$ the first assertion follows.

The second one is clear since by Theorem \ref{theo04} cohomologies in degree $1$ agree. Hence decomposable cohomologies in degrees $2$ are isomorphic. Consequently, also the relations imposed on $\Lambda V_0$ in degree $2$ have to agree.
\end{prf}

\section{Proof of Theorem \ref{theoB}}\label{sec02}

We may transcribe the Lie algebras in the language of minimal Sullivan algebras. We prove Theorem \ref{theoB} case by case.
\begin{theo}\label{theo01}
For $k\geq 1$ the algebra $(\Lambda V,\dif)$ with
\begin{align*}
V_0&=\langle x_1, \ldots , x_{2k}\rangle,\\
V_1&=\langle n_1, \ldots, n_{2k-1}\rangle,\\
V_2&=\langle m\rangle
\end{align*}
and $\dif x_i=0$ for $1\leq i\leq 2k$, $\dif n_i=x_{i}x_{i+1}$ and $\dif m=\sum_{i=1}^{2k-1} x_{i}n_i$ is graded and rigid.
\end{theo}
\begin{prf}
It is clear that the corresponding nilpotent Lie algebra is graded due to the grading on the differential, i.e.~since the differential reduces lower degree by exactly one.

\vspace{5mm}

\step{1}
We consider a second Lie algebra and its associated graded Lie algebra. We encode both by their minimal models $(\Lambda W,\dif)$ and $(\Lambda W,\bar\dif)$. We show that if $(\Lambda W,\bar\dif)$ is isomorphic to $(\Lambda V,\bar \dif)=(\Lambda V, \dif)$, then $W=W_0\oplus W_1 \oplus W_2$, i.e.~$W$ is $3$-stage, and that $W_i=V_i$ for $0\leq i\leq 2$.

The observation is easily seen as follows. First we identify $(\Lambda W,\bar \dif)$ with $(\Lambda V,\dif)$. By the construction of the minimal model $(\Lambda W,\bar \dif)$ out of $(\Lambda W,\dif)$, i.e.~the fact that $\dif-\bar \dif$ decreases lower degree by two, we see (using the identification above) that
\begin{align*}
\dif W_0&=0 \\
\dif n_i&=  x_{i}x_{i+1}\\
\dif m&=\sum_{i=1}^{2k-1} x_{i}n_i + p
\end{align*}
where $p\in \Lambda^2 \langle x_i\rangle$. This yields the detected form of $(\Lambda W,\dif)$.

\vspace{5mm}

\step{2} It remains to show that $(\Lambda V,\dif)\cong (\Lambda W,\dif)$ with the differential of the latter as specified above.

We know that $(\Lambda (V_0\oplus V_1),\dif)$ and $(\Lambda (W_0\oplus W_1),\dif)$ are isomorphic by Lemma \ref{lemma02}.
It follows that $(\Lambda W,\dif)$ splits as a real fibration
\begin{align*}
(\Lambda (V_0\oplus V_1),\dif)\hto{} (\Lambda W,\dif)\to (\Lambda \langle m\rangle,0)
\end{align*}
In other words, $(\Lambda W,\dif)\cong (\Lambda (V_0\oplus V_1\oplus \langle m\rangle),\dif)$ is completely determined by $(\Lambda (V_0\oplus V_1),\dif)$ and the differential on $m$ given by $\dif m\in \Lambda (V_0\oplus V_1)$. Additionally, we know that $\bar \dif m=\sum_{i=1}^{2k-1} x_{i}n_i$.

Let us now see that up to isomorphism $\bar \dif(m)=\dif(m)$ which proves the result. Again consider the general form of $\dif (m)$ as given by
\begin{align*}
\dif m&=\sum_{i=1}^{2k-1} x_{i}n_i + p
\end{align*}
with $p\in \Lambda^2 \langle x_i\rangle$.

Now write $p$ in the form
\begin{align*}
p=&x_1(t_{1,2}x_2+t_{1,3}x_3+\ldots +t_{1,2k} x_{2k})\\
+&x_2(t_{2,3}x_3 +\ldots + t_{2,2k} x_{2k})\\
&\vdots\\
+&x_{2k-1}(t_{2k-1,2k}x_{2k})
\end{align*}
with coefficients $t_{i,j}\in \rr$.

This permits us to finally specify an isomorphism $f\co (\Lambda W,\dif)\to (\Lambda V,\dif)$ by
\begin{align*}
f(x_i)&=x_i\\
f(n_j)&=n_j-(t_{j,j+1} x_{j+1} + t_{j,j+2} x_{j+2} + \ldots + t_{j,2k} x_{2k})\\
f(m)&=m
\end{align*}
for all $i,j$. This is clearly an isomorphism once it commutes with differentials. Since $\dif x_j=0$, it remains only to check that $f(\dif m)=\dif f(m)$. We compute
\begin{align*}
f(\dif m)=&f\bigg(\sum_{i=1}^{2k-1} x_{i}n_i + p\bigg)\\
=&f\bigg(\bigg(\sum_{i=1}^{2k-1} x_{i}n_i\bigg)+(x_1(t_{1,2}x_2+t_{1,3}x_3+\ldots \\&+t_{1,2k} x_{2k})+ \ldots +x_{2k-1}(t_{2k-1,2k}x_{2k}))\bigg)
\\=&\big( x_1(n_1-(t_{1,2} x_{2} + t_{1,3} x_{3} + \ldots + t_{1,2k} x_{2k}))\\
&+ x_2(n_2-(t_{2,3} x_{2} + t_{2,4} x_{4} + \ldots + t_{2,2k} x_{2k}))\\
&\vdots\\
&+x_{2k-1}(n_{2k-1}-t_{2k-1,2k}x_{2k})\big)
\\&+\big((x_1(t_{1,2}x_2+t_{1,3}x_3+\ldots +t_{1,2k} x_{2k})+ \ldots +x_{2k-1}(t_{2k-1,2k}x_{2k}))\big)
\\=&\sum_{i=1}^{2k-1} x_in_i
\\=&\dif(f(m))
\end{align*}
\end{prf}

Now we present a non-graded rigid example.  In the proof we shall not only use Theorem \ref{theo04} as in the proof of Theorem \ref{theo01}, but also shall we draw on Theorem \ref{theo05}.

We conjecture that the following theorem is true for all $k\geq 2$.
\begin{theo}\label{theo02}
For $k\in \{2,3\}$ the algebra $(\Lambda V,\dif)$ with
\begin{align*}
V_0&=\langle x_1, \ldots , x_{2k}, x_{2k+1}\rangle,\\
V_1&=\langle n_1, \ldots, n_{2k-1}\rangle,\\
V_2&=\langle m\rangle
\end{align*}
and $\dif x_i=0$ for $1\leq i\leq 2k+1$, $\dif n_i=x_{i}x_{i+1}$ and
\begin{align*}
\dif m=\bigg(\sum_{i=1}^{2k-1} x_{i}n_i\bigg)+x_{2k}\cdot x_{2k+1}
\end{align*}
is non-graded and rigid.
\end{theo}
\begin{prf}
Again by construction of the model of the graded algebra it is clear that the algebra is not graded, and no isomorphic image of it is either---this is basically due to the fact that $x_{2k}x_{2k+1}$ is not exact. More precisely, it is a direct check that no such isomorphism to its associated graded algebra exists. Such an isomorphism would have to alter the $n_i$ by summands of $x_i$ such that $x_{2k}x_{2k+1}$ lies in the image of $\sum_{i=1}^{2k-1}x_in_i$. Since every summand in this sum has a coefficient $x_i$ with $i<2k$, this is impossible. In the last step of the proof we actually show, that both algebras even do not have isomorphic cohomologies for $k\in \{2,3\}$.

The associated graded algebra has identical minimal model $(\Lambda V,\bar \dif)$, \linebreak[4]i.e.~$\dif(x_i)=\bar \dif(x_i)$, $\dif(n_j)=\bar \dif(n_j)$, except for $\bar\dif(m)=\dif(m)-x_{2k}\cdot x_{2k+1}$. As in the proof of Theorem \ref{theo01} we proceed in two steps. First we show that any other quasi-isometric nilpotent Lie algebra encoded by its minimal model $(\Lambda W,\dif)$ is $3$-stage; second we show that it has to be isomorphic. One major difference in this second step will be that we need to use Theorem \ref{theo04} and not only Theorem \ref{theo05}.

\vspace{5mm}

\step{1} We proceed as in Step 1 of Theorem \ref{theo01}. We consider $(\Lambda W,\dif)$ with $(\Lambda W,\bar \dif)=(\Lambda V,\bar \dif)$. Again, it follows that
\begin{align*}
\dif|_{W_0\oplus W_1}&= \bar \dif_{W_0\oplus W_1}\\
\dif m&=\bigg(\sum_{i=1}^{2k-1} x_{i}n_i\bigg) + p
\end{align*}
with $p\in \Lambda^2 \langle x_i\rangle$.

\vspace{5mm}

\step{2} Now let us see that $(\Lambda W,\dif)$ and $(\Lambda V,\dif)$ are isomorphic. For reasons of lower degree $(\Lambda (W_0\oplus W_1),\dif)=(\Lambda (W_0\oplus W_1),\bar \dif)$, which by Lemma \ref{lemma02} is isomorphic to $(\Lambda (V_0\oplus V_1),\dif)$. Again, as in Theorem \ref{theo02}, we derive that $(\Lambda W,\dif)$ is completely determined by $(\Lambda (V_0\oplus V_1),\dif)$ together with the differential $\dif m\in \Lambda (V_0\oplus V_1)$. Moreover, we know that $\bar \dif m=\sum_{i=1}^{2k-1} x_{i}n_i$.

Thus it remains to discuss $p$. Write $p$ in the form
\begin{align*}
p=&x_1(t_{1,2}x_2+t_{1,3}x_3+\ldots +t_{1,2k+1} x_{2k+1})\\
+&x_2(t_{2,3}x_3 +\ldots + t_{2,2k+1} x_{2k+1})\\
&\vdots\\
+&x_{2k-1}(t_{2k-1,2k}x_{2k}+t_{2k-1,2k+1}x_{2k+1})\\
+&x_{2k}(t_{2k,2k+1}x_{2k+1})
\end{align*}
with coefficients $t_{i,j}\in \rr$.

We use this to simplify our algebra. That is, let us define the algebra $(\Lambda W,\tilde \dif)$ where the differential is given as
\begin{align*}
\tilde \dif(x_i)&=0\\
\tilde \dif(n_j)&=\dif(n_j)=x_jx_{j+1}\\
\tilde \dif(m)&=\bar \dif(m)+ t_{2k,2k+1} x_{2k}x_{2k+1}=\sum_{i=1}^{2k-1} x_{i}n_i + t_{2k,2k+1} x_{2k}x_{2k+1}
\end{align*}
for all $i,j$.

We specify an isomorphism $f\co (\Lambda W,\dif)\to (\Lambda W,\tilde\dif)$ by
\begin{align*}
f(x_i)&=x_i\\
f(n_j)&=n_j-(t_{j,j+1} x_{j+1} + t_{j,j+2} x_{j+2} + \ldots + t_{j,2k+1} x_{2k+1})\\
f(m)&=m
\end{align*}
for all $i,j$. Again, this is clearly an isomorphism once it commutes with differentials. Since $\dif x_j=0$, it remains only to check that $f(\dif m)=\tilde\dif f(m)$. We compute
\begin{align*}
f(\dif m)=&f\bigg(\sum_{i=1}^{2k-1} x_{i}n_i + w\bigg)\\
=&f\bigg(\bigg(\sum_{i=1}^{2k-1} x_{i}n_i\bigg)+(x_1(t_{1,2}x_2+t_{1,3}x_3+\ldots \\&+t_{1,2k+1} x_{2k+1})+ \ldots +x_{2k-1}(t_{2k-1,2k}x_{2k}  +t_{2k-1,2k+1}x_{2k+1})\\&+x_{2k}(t_{2k,2k+1}x_{2k+1}))\bigg)
\\=&\big( x_1(n_1-(t_{1,2} x_{2} + t_{1,3} x_{3} + \ldots + t_{1,2k+1} x_{2k+1}))\\
&+ x_2(n_2-(t_{2,3} x_{2} + t_{2,4} x_{4} + \ldots + t_{2,2k+1} x_{2k+1}))\\
&\vdots\\
&+x_{2k-1}(n_{2k-1}-t_{2k-1,2k}x_{2k}-t_{2k-1,2k+1}x_{2k+1})\big)
\\&+\big((x_1(t_{1,2}x_2+t_{1,3}x_3+\ldots +t_{1,2k+1} x_{2k+1})+ \ldots \\&+x_{2k-1}(t_{2k-1,2k}x_{2k}+t_{2k-1,2k+1}x_{2k+1})+x_{2k}(t_{2k,2k+1}x_{2k+1}))\big)
\\=&\sum_{i=1}^{2k-1} x_in_i +t_{2k,2k+1}x_{2k}x_{2k+1}
\\=&\tilde\dif(f(m))
\end{align*}
Set $t:=t_{2k,2k+1}$. For the rest of this proof it remains to distinguish two cases: First $t\neq 0$, second $t=0$.

\case{1} If $t\neq 0$, we apply the morphism $g\co (\Lambda W,\tilde \dif)\to (\Lambda V,\dif)$ induced by $x_{2k+1}\mapsto (1/t)\cdot x_{2k+1}$. It is easily seen to commute with differentials---actually only $\dif(g(m))=g(\tilde\dif(m))$ need to be checked and we have
\begin{align*}
\dif(g(m))=&\sum_{i=1}^{2k-1} x_in_i +x_{2k}x_{2k+1}\\
=&g\bigg(\sum_{i=1}^{2k-1} x_in_i +tx_{2k}x_{2k+1}\bigg)\\
=&g(\tilde \dif (m))
\end{align*}
Thus this morphism is actually an isomorphism of minimal Sullivan algebras.

\vspace{3mm}

\case{2} If $t=0$, we see that $(\Lambda W,\dif)$ is isomorphic to $(\Lambda W,\tilde \dif)$ which  is identical to $(\Lambda V,\bar \dif)$. In this case we obtain a contradiction from Theorem \ref{theo04} once we have shown that the algebras $H(\Lambda V,\dif)$ and $H(\Lambda V,\bar \dif)$ for $k\in \{2,3\}$ are not graded isomorphic.

We do so by observing that the number of cohomology algebra generators in degree $3$ in all these cases differs by one. For this we compute the primitive cohomology using the computer algebra system \textsc{Sage-7.6}.

For the convenience of the reader we describe the basis of primitive cohomology in degree $3$ in these cases:
They are
\begin{align*}
[x_4\cdot x_5\cdot n_2 + x_1\cdot n_1\cdot n_2 - x_2\cdot x_3\cdot m],\\
[x_4\cdot x_5\cdot n_1 - x_3\cdot n_1\cdot n_3 - x_1\cdot x_2\cdot m],\\
[x_3\cdot n_2\cdot n_3],\\
[x_3\cdot n_1\cdot n_2 + x_1\cdot x_5\cdot n_3 - x_1\cdot x_3\cdot m],\\
[x_2\cdot n_1\cdot n_2],\\
[x_1\cdot n_1\cdot n_2 + x_2\cdot x_5\cdot n_3 - x_2\cdot x_3\cdot m]]
\end{align*}
for $\dif$ and $k=2$,
\begin{align*}
[x_3\cdot n_2\cdot n_3],\\
[x_3\cdot n_1\cdot n_3 + x_1\cdot x_2\cdot m],\\
[x_3\cdot n_1\cdot n_2 - x_1\cdot x_3\cdot m],\\
[x_2\cdot n_1\cdot n_2],\\
[-x_1\cdot n_1\cdot n_2 + x_2\cdot x_3\cdot m]
\end{align*}
for $\bar \dif$ and $k=2$,
\begin{align*}
[x_6\cdot x_7\cdot n_4 + x_1\cdot n_1\cdot n_4 + x_2\cdot n_2\cdot n_4 + x_3\cdot n_3\cdot n_4 - x_4\cdot x_5\cdot m],\\
[x_6\cdot x_7\cdot n_3 + x_1\cdot n_1\cdot n_3 + x_2\cdot n_2\cdot n_3 - x_5\cdot n_3\cdot n_5 - x_3\cdot x_4\cdot m],\\
[x_6\cdot x_7\cdot n_2 + x_1\cdot n_1\cdot n_2 - x_4\cdot n_2\cdot n_4 - x_5\cdot n_2\cdot n_5 - x_2\cdot x_3\cdot m],\\
[x_6\cdot x_7\cdot n_1 - x_3\cdot n_1\cdot n_3 - x_4\cdot n_1\cdot n_4 - x_5\cdot n_1\cdot n_5 - x_1\cdot x_2\cdot m],\\
[x_5\cdot n_4\cdot n_5],\\
[x_4\cdot n_3\cdot n_4],\\
[x_1\cdot n_1\cdot n_4 + x_2\cdot n_2\cdot n_4 + x_3\cdot n_3\cdot n_4 + x_4\cdot x_7\cdot n_5 - x_4\cdot x_5\cdot m],\\
[x_3\cdot n_2\cdot n_3],\\
[x_2\cdot n_1\cdot n_2]
\end{align*}
for $\dif$ and $k=3$,
\begin{align*}
[x_5\cdot n_4\cdot n_5],\\
[ -x_1\cdot n_1\cdot n_3 - x_2\cdot n_2\cdot n_3 + x_5\cdot n_3\cdot n_5 + x_3\cdot x_4\cdot m],\\
[  -x_1\cdot n_1\cdot n_2 + x_4\cdot n_2\cdot n_4 + x_5\cdot n_2\cdot n_5 + x_2\cdot x_3\cdot m],\\
[  x_3\cdot n_1\cdot n_3 + x_4\cdot n_1\cdot n_4 + x_5\cdot n_1\cdot n_5 + x_1\cdot x_2\cdot m],\\
[  x_4\cdot n_3\cdot n_4],\\
[  -x_1\cdot n_1\cdot n_4 - x_2\cdot n_2\cdot n_4 - x_3\cdot n_3\cdot n_4 + x_4\cdot x_5\cdot m],\\
[  x_3\cdot n_2\cdot n_3],\\
  [x_2\cdot n_1\cdot n_2]
\end{align*}
for $\bar \dif$ and $k=3$,

\vspace{5mm}

The consequence of the Case 1, Case 2 distinction is that the algebras $(\Lambda V,\dif)$ and $(\Lambda W,\dif)$ are isomorphic; as claimed in the assertion.
\end{prf}

\vspace{5mm}

Denote by $(\Lambda~_{k}^{l}F,\dif):=(\Lambda(^{l}F_0\oplus \ldots \oplus~ ^{l}F_k ),\dif)$ the minimal model of the free nilpotent Lie algebra over $l$ generators in lower degree $0$; it is supposed to be ``free nilpotent till lower degree $k$''. By this we mean that $H_{\leq k-1}(\Lambda~_{k}^{l}F)=0$, and first non-trivial cohomology appears in lower degree $k$.
\begin{theo}\label{theo03}
Let $V\In ^lF_{k+1}$ be a subspace.
The Lie algebra $(\Lambda~(_{k}^{l}F\oplus V),\dif)$ with the differential on $V$ defined by restriction of the differential on $^lF_{k+1}$ is graded and rigid.
\end{theo}
\begin{prf}

Let $(\Lambda W,\dif)$ be a minimal model such that the associated graded models $(\Lambda W,\bar\dif)\cong (\Lambda~(_{k}^{l}F\oplus V),\bar \dif)$ are isomorphic and such that $H(\Lambda W,\dif)\cong H(\Lambda~(_{k}^{l}F\oplus V),\dif)$. We observe that $ (\Lambda~(_{k}^{l}F\oplus V),\dif)= (\Lambda~(_{k}^{l}F\oplus V),\bar \dif)$ is Carnot-graded. It follows that
\begin{align*}
(\Lambda~(_{k}^{l}F\oplus V),\dif)\cong (\Lambda W,\dif)
\end{align*}
up to perturbations of $\dif$ on $W$, i.e.~up to summands in $\dif(w)$ for $w\in W$ of lower degree smaller than the lower degree of $\bar \dif w$. Up to isomorphism we can assume that this is a perturbation by closed non-exact classes. Indeed, we proceed inductively of lower degree. Let $z\in W$ be of minimal lower degree such that $\dif(z)=z'+z''$ is not lower graded homogeneous with $z'$ the homogeneous and $z''$ the non-homogeneous part. Then
\begin{align*}
0=\dif^2(z)=\bar\dif(z'+z'')=\bar\dif^2 z' +\bar \dif z''= \bar \dif z''=\dif z''
\end{align*}
Up to isomorphism we may assume it to be non-exact. Then however, we observe that three are no non-trivial cohomology classes by which we could perturb the differential, since cohomology vanishes in lower degrees at most $k$ and ordinary degree $2$. We may iterate this argument to see that no perturbation at all can be realised. Consequently, the isomorphism $(\Lambda~(_{k}^{l}F\oplus V),\dif)\cong (\Lambda W,\dif)$ holds in general, and the associated Lie algebra is rigid.
\end{prf}

\vspace{5mm}

Now we prove the last case of our main Theorem \ref{theoB}, yielding another non-graded example.

\begin{theo}\label{theo04}
The algebra $(\Lambda V,\dif)$ with
\begin{align*}
V_0&=\langle x_1,x_2,x_3,x_4\rangle,\\
V_1&=\langle n_1,n_2,n_3,n_4,n_5\rangle,\\
V_2&=\langle m\rangle
\end{align*}
and $\dif x_i=0$ for $1\leq i\leq 4$, $\dif n_1=x_1x_2, \dif n_2=x_1x_3, \dif n_3=x_2x_3, \dif n_4=x_1x_4, n_5=x_2x_4$ and
\begin{align*}
\dif m=x_1n_1+x_1n_2+x_2n_3+x_3x_4
\end{align*}
is non-graded and rigid.
\end{theo}
\begin{prf}
We proceed as in the previous proofs using analog arguments. Thus, again by construction of the model of the graded algebra it is clear that the algebra is not graded.
Let $(\Lambda W,\dif)$ an algebra with isomorphic graded algebra and isomorphic cohomology. Then $(\Lambda W,\bar \dif)$ equals $(\Lambda V,\dif)$ except for $\dif m=x_1n_1+x_1n_2+x_2n_3$. Since we have that
\begin{align*}
H^{\geq 2}_{0}(\Lambda W,\bar\dif)= \langle x_1x_4\rangle
\end{align*}
we deduce that the only perturbation of $\dif m$ in $(\Lambda W,\dif)$ up to isomorphism is of the form $\dif m=x_1n_1+x_1n_2+x_2n_3+t \cdot x_1x_4$.

Again we have to differ two cases. If $t\neq 0$ we consider the isomorphism of minimal models induced by $x_4\mapsto 1/t \cdot x_4$, $n_4\mapsto 1/t \cdot n_4$, $n_5\mapsto 1/t \cdot n_5$, $m\mapsto m$, which allows us to assume that $t=1$ and $(\Lambda W,\dif)\cong (\Lambda V,\dif)$.

If $t=0$ we compute using \textsc{Sage-7.6} that $H(\Lambda V,\dif)$ has $31$ algebra generators in degree $3$ whereas $H(\Lambda V,\bar \dif)$ has $27$ of them in degree $3$. Thus both cohomology algebras are not isomorphic.

Consequently, either $(\Lambda V,\dif)$ and $(\Lambda W,\dif)$ (under the assumption of isomorphic graded algebras) are isomorphic or have distinct cohomology algebras already. This proves the assertion.
\end{prf}

\vspace{5mm}

We leave it to the interested reader to vary these construction principles and to find new series of rigid examples.


\section{A non-example}\label{sec03}

For completeness we provide an example two nilpotent Lie algebras which have isomorphic associated graded Lie algebras and isomorphic cohomology, but which are not isomorphic themselves. This example seems to be well-known already, see

{\fontsize{6}{9}\selectfont
\centering
\begin{verbatim}
https://mathoverflow.net/questions/122238/malcev-lie-algebra-and-associated-graded-lie-algebra/122267#122267
\end{verbatim}
}

This shows that the combination of the results by Sauer and Pansu does not uniquely characterise the nilpotent Lie algebra, and, at least a priori, is not strong enough to completely solve the main conjecture.

For this we define the following two $5$-dimensional $4$-stage minimal models $(\Lambda V,\dif_{1/2})=(\Lambda V^1,\dif_{1/2})$ with $V^1=\langle a_1,a_2,b,c,d\rangle$ with
\begin{align*}
\dif_1 a_1=\dif_1 a_2&=0 \\
\dif_1 b&=a_1a_2\\
\dif_1 c&=a_1b\\
\dif_1 d&=a_1c\\
\end{align*}
respectively
\begin{align*}
\dif_2 a_1=\dif_2 a_2&=0 \\
\dif_2 b&=a_1a_2\\
\dif_2 c&=a_1b\\
\dif_2 d&=a_1c+ a_2b\\
\end{align*}

The associated graded algebra is just the same in each case given by $(\Lambda V,\dif_1)$ itself.
It is easy to see that the algebras are not isomorphic. By degree reasons any such isomorphic comes from an automorphism of $V$. However, it is obvious that we cannot specify an image of $a_1$ and of $c$ such that $a_1c\mapsto a_1c+a_2d$.

Indeed for
\begin{align*}
a_1\mapsto& k_1a_1+k_2a_2+k_3b+k_4c+k_5d\\
c\mapsto &l_1a_1+l_2a_2+l_3b+l_4c+l_5d
\end{align*}
with real $k_i,l_i$ we compute
\begin{align*}
&a_1c\mapsto \\&(k_1l_2-k_2l_1)a_1a_2 + (k_1l_3-k_3l_1)a_1b + (k_1l_4-k_4l_1)a_1c+ (k_1l_5-k_5l_1)a_1d + \\
&(k_2l_3-k_3l_2)a_2b + (k_2l_4-k_4l_2)a_2c + (k_2l_5-k_5l_2)a_2d+\\
&(k_3l_4-k_4l_3)bc + (k_3l_5-k_5l_3)bd +\\
&(k_4l_5-k_5l_4)cd
\end{align*}
Thus the following terms need to vanish:
\begin{align*}
&(k_1l_2-k_2l_1), (k_1l_3-k_3l_1), (k_1l_5-k_5l_1) \\
&(k_2l_3-k_3l_2), (k_2l_4-k_4l_2), \\
&(k_3l_4-k_4l_3), (k_3l_5-k_5l_3),\\
&(k_4l_5-k_5l_4)
\end{align*}
and $(k_1l_4-k_4l_1)=(k_2l_5-k_5l_2)=1$. We may assume that $k_1=l_4=1$ so that $k_4=0$ or $l_1=0$. It follows that $k_4=0\Rightarrow k_2=k_3=k_4=k_5=0 \vee l_4=0$ and $l_1=0 \Rightarrow k_1=0 \vee l_1=l_2=l_3=l_5=0$. That is, we have that $(k_1=1,k_2=k_3=k_4=k_5=0)$ and $l_1=l_2=l_3=l_5=0, l_4=1$. Then, however, $(k_2l_5-k_5l_2=0$; a contradiction.

\vspace{5mm}

Let us now describe the cohomology ring in both cases and see that it is isomorphic. We denote generators of the cohomology algebra together with the morphism upon them.
\begin{align*}
[a_1] \mapsto & [a_1]\\
[a_2] \mapsto & [a_2]\\[1em]
[a_2b]\mapsto & [a_2b] ~~(=[a_1c])\\
[bc-a_2d]\mapsto &[bc-a_2d]\\[1em]
[a_1bc] \mapsto & [a_1bc]\\
[a_1cd]\mapsto & [a_1cd-a_2bd]\\[1em]
[a_1bcd]\mapsto &[a_1bcd]
\end{align*}
The morphism is obviously an isomorphism of modules. It is easy but tedious that this morphism is actually compatible with the multiplicative structure. Thus it is an isomorphism of cohomology algebras.

It is now easy to describe the nilpotent Lie algebras, these are the models of, namely the algebras on the vector space $\mathfrak{g}:=\langle a_1,a_2,b,c,d\rangle$ with non-trivial Lie brackets once given by
\begin{align*}
[a_1,a_2]=-b, [a_1,b]=-c, [a_1,c]=-d
\intertext{and second by}
[a_1,a_2]=-b, [a_1,b]=-c, [a_1,c]=-d, [a_2,b]=-d
\end{align*}

\vspace{5mm}

Hence we have proved
\begin{theo}\label{theo07}
There are examples of nilpotent Lie algebras which have isomorphic graded algebras and isomorphic cohomology, but which are not isomorphic themselves.
\end{theo}



\def\cprime{$'$}


\vfill

\begin{center}
\noindent
\begin{minipage}{\linewidth}
\small \noindent \textsc
{Manuel Amann} \\
\textsc{Institut f\"ur Mathematik}\\
\textsc{Differentialgeometrie}\\
\textsc{Universit\"at Augsburg}\\
\textsc{Universit\"atsstra\ss{}e 14 }\\
\textsc{86159 Augsburg}\\
\textsc{Germany}\\
[1ex]
\textsf{manuel.amann@math.uni-augsburg.de}\\
\textsf{www.math.uni-augsburg.de/prof/diff/arbeitsgruppe/amann/}
\end{minipage}
\end{center}

\end{document}